\documentclass[11pt,reqno]
{amsart} 
\usepackage{verbatim}
\usepackage{amssymb}
\usepackage{enumerate}
\usepackage[active]{srcltx}

\newif\ifdeveloping

\def\myheads#1;#2;{
\pagestyle{myheadings} \markboth{{\sc\hfill
#1\hfill\protect\makebox[0cm][r]{\rm\today}}}
{{\sc\protect\makebox[0cm][l]{\rm\today}\hfill #2\hfill}} }

\newcommand{\setm}{\setminus}
\newcommand{\empt}{\emptyset}
\newcommand{\subs}{\subset}

\def\<{\left\langle}
\def\>{\right\rangle}

\def\OO{{\omega}}
\def\br#1;#2;{\bigl[ {#1} \bigr]^ {#2} }

\def\ooseq#1;#2;{\< {#1}_{#2}:{#2}<\oo\>}
\def\ooset#1;#2;{\{ {#1}_{#2}:{#2}<\oo\}}
\def\seq#1;#2;#3;{\< {#1}_{#2}:{#2}<#3\>}
\def\set#1;#2;#3;{\{ {#1}_{#2}:{#2}<#3\}}
\def\oseq#1;#2;{\< {#1}_{#2}:{#2}<\OO\>}
\def\oset#1;#2;{\{ {#1}_{#2}:{#2}<\OO\}}
\def\oosequ#1;#2;{\< {#1}^{#2}:{#2}<\oo\>}
\def\oosetu#1;#2;{\{ {#1}^{#2}:{#2}<\oo\}}
\def\sequ#1;#2;#3;{\< {#1}^{#2}:{#2}<#3\>}
\def\setu#1;#2;#3;{\{ {#1}^{#2}:{#2}<#3\}}
\def\osequ#1;#2;{\< {#1}^{#2}:{#2}<\OO\>}
\def\osetu#1;#2;{\{ {#1}^{#2}:{#2}<\OO\}}

\newcommand{\al}{\alpha}
\newcommand{\be}{\beta}

\ifdeveloping
\usepackage[notref,notcite]{showkeys}
\fi
\usepackage{enumerate}

\newcommand{\prtime}{{\count0=\time\divide\count0 by 60
\count1=-\count0\multiply\count1 by 60 \advance\count1 by \time
\the\count0:\the\count1} }

\def\myheads#1;#2;{
\pagestyle{myheadings} \markboth{{\sc\hfill
#1\hfill\protect\makebox[0cm][r]{\rm\today; \prtime}}}
{{\sc\protect\makebox[0cm][l]{\rm\today;\ \prtime}\hfill
#2\hfill}} \thispagestyle{myheadings} }

\newcommand{\de}{\delta}
\newcommand{\ka}{\kappa}
\newcommand{\la}{\lambda}

\newtheorem{theorem}{Theorem}[section]
\newtheorem{proposition}[theorem]{Proposition}
\newtheorem{definition}[theorem]{Definition}
\newtheorem{definitions}[theorem]{Definitions}

\newtheorem{corollary}[theorem]{Corollary}

\newcommand{\ga}{\gamma}
\newcommand{\om}{\omega}
\def\<{\left\langle}
\def\>{\right\rangle}

\title{ Constructions of Lindel\"{o}f scattered P-spaces}

\author[J. C. Martinez]{Juan Carlos Mart\'{\i}nez}
\address{Facultat de Matem\`atiques  i Inform\`atica\\ Universitat de Barcelona \\ Gran
 Via 585 \\ 08007 Barcelona, Spain}
\email{jcmartinez@ub.edu}
\author[L. Soukup]{
Lajos Soukup }
\address{Alfr{\'e}d R{\'e}nyi Institute of Mathematics }
\email{soukup@renyi.hu}

\thanks{The preparation of this paper  was  supported by Hungarian National Foundation for Scientific Research grant no. 129211 and 
 Spanish MICINN
Grant PID2020-116773GB-100}

\begin{document}

\footnotetext[1] { 2010 {\em Mathematics Subject Classification}.  54A25, 54A35, 54G12, 03E35
\\ \hspace*{5mm} {\em Keywords and phrases}. Lindel\"{o}f scattered P-space, Cantor-Bendixson height, Cantor-Bendixson width .}

\begin{abstract} 
  We construct  locally Lindelöf scattered P-spaces 
  (LLSP spaces, in short)
  with prescribed widths and heights under different set-theoretic assumptions.

  We prove that there is an LLSP space of width $\om_1$ and height $\om_2$ and 
  that it is relatively consistent with ZFC that there is an LLSP space 
  of width $\om_1$ and height $\om_3$. 
  Also, we prove a stepping up theorem that, for every cardinal $\la \geq \om_2$, 
  permits us to construct from an LLSP space
    of width $\om_1$ and height $\la$ satisfying certain additional properties an
 LLSP space of width $\om_1$ and height $\al$ for every ordinal $\al < \la^+$. 
    Then, we obtain as consequences of the above results the following theorems:


\vspace{1mm} (1) For every ordinal $\al < \om_3$ there is an LLSP space  
 of width $\om_1$ and height $\al$.

\vspace{1mm} (2) It is relatively consistent with ZFC that 
there is an LLSP space of width $\om_1$ and height $\al$ for every ordinal $\al < \om_4$.

\end{abstract}

\maketitle

\section{Introduction}

The cardinal sequence of a scattered space  is the sequence of the cardinalities of its Cantor-Bendixson levels.
The investigation of the cardinal sequences of different classes of topological spaces is a classical problem of set theoretic topology. 
Many  important results were proved in connection with the cardinal sequences of 
locally compact scattered  (LCS, in short) spaces, see e.g.  
\cite{ba2007,Ba2002,bs87,ErVe2010,jw78,Juwe2006,Ju85,Lag77,Ma1992,ma01,ma03,R76a,Ro2002,Ro85}.  
In \cite{JSSSh2004} a complete characterization of the cardinal sequences of the 
0-dimensional, of the  regular, and of the Hausdorff spaces was given.  

Recall that a topological space $X$ is a {\em P-space}, if the intersection of every countable family of open sets in $X$ is open in $X$. 
The aim of this paper is to start the systematic investigation of cardinal sequences 
of locally Lindelöf scattered P-spaces. We will see that several methods applied to LCS spaces 
can be applied here, but typically we should face more serious technical problems.

\vspace{2mm}  If $X$ is a topological space and $\al$ is an ordinal, we denote by $X^{\al}$ the $\al$-th Cantor-Bendixson derivative of $X$. Then, $X$ is {\em scattered} if $X^{\al} = \emptyset$ for some ordinal $\al$. Assume that $X$ is a scattered space. We define the {\em  height } of  $X$ by

$$\mbox{ht}(X) = \mbox { the least ordinal } \al \mbox{ such that } X^{\al} = \emptyset.$$

   \vspace{1mm}
   \noindent For $\al < \mbox{ht}(X)$, we write $I_{\al}(X) = X^{\al}\setminus X^{\al + 1}$. If $x\in I_{\al}(X)$, we say that $\al$ is the {\em level} of $x$ and we write $\rho(x,X) = \al$, or simply $\rho(x) = \al$ if no confusion can occur. Note that $\rho(x) = \al$ means that $x$ is an accumulation point of $I_{\beta}(X)$ for $\beta < \al$ but $x$ is not an accumulation point of $X^{\al} = \bigcup \{I_{\be}(X) : \be \geq \al \}$. We define the {\em width} of $X$ as

   $$\mbox{wd}(X) = \mbox{sup}\{ |I_{\al}(X)| : \al < \mbox{ht}(X) \}.$$

   If $X$ is a scattered space, $x\in X$ and $U$ is a neighbourhood of $x$, we say that $U$ is a {\em cone on} $x$, if $x$ is the only point in $U$ of level $\geq \rho(x,X)$.

\vspace{1mm}
   By an {\em LLSP space}, we mean a locally Lindel\"{o}f, scattered, Hausdorff P-space.

   \begin{proposition}\label{pr:mi}
    An LLSP space  is 0-dimensional.    
    \end{proposition}

  \begin{proof}
    By \cite[Proposition 4.2(b)]{Mi}, a Lindelöf Hausdorff P-space $X$  is normal, so a locally Lindelöf Hausdorff P-space is regular.  
    Thus, by  \cite[Corollary 3.3]{Mi}, $X$ is 0-dimensional.
    \end{proof}

So, by  Proposition  \ref{pr:mi} above,  
   if  $X$ is an LLSP space, $x\in X$ and ${\mathbb B}_x$ is a neighbourhood basis of $x$, we may assume that every $U\in {\mathbb B}_x$ is a Lindel\"{o}f clopen cone on $x$.

 \bigskip

 It was proved by Juh\'asz and Weiss in \cite{jw78} that for every ordinal ${\alpha} < {\omega}_2$
there is an LCS space of height ${\alpha}$ and width ${\omega}$. Then, we will transfer this
theorem to the setting of LLSP spaces, showing that for every ordinal ${\alpha} < {\omega}_3$
there is an LLSP space of height ${\alpha}$ and width ${\omega}_1$.

To obtain an LCS space of height ${\omega}_1$ and width ${\omega}$,
 in \cite{jw78} Juhász and Weiss, using transfinite recursion, constructed a sequence $\<X_{\alpha}:{\alpha}\le{\omega}_1\>$ of LCS spaces such that 
$X_{\alpha}$ had height ${\alpha}$ and width ${\omega}$, and for ${\alpha}<{\beta}$, the space
$X_{\alpha}$ was just the first ${\alpha}$ Cantor-Bendixson levels of $X_{\beta}$.

Since $X_{\alpha}$ is dense in $X_{{\alpha}+1}$, Juhász and Weiss had to guarantee
that $X_{\alpha}$ is not compact. But it was automatic, because if ${\alpha}=\gamma+1$, then 
$X_{\alpha}$ had a top  infinite Cantor-Bendixson level, so $X_{\alpha}$  was not compact.
If ${\alpha}$ is a limit ordinal, then the open cover $\{X_{\xi}:{\xi}<{\alpha}\}$ witnessed
that $X_{\alpha}$ is not compact. 

What happens if we try to adopt that approach for LLSP spaces? 
To obtain an LLSP space of height ${\omega}_2$ and width ${\omega}_1$,
we can try, using transfinite recursion, to construct a sequence 
$\<X_{\alpha}:{\alpha}\le{\omega}_2\>$ of LLSP spaces such that 
$X_{\alpha}$ has height ${\alpha}$ and width ${\omega}_1$, and for ${\alpha}<{\beta}$, the space
$X_{\alpha}$ is just the first ${\alpha}$ levels of $X_{\beta}$.

Since $X_{\alpha}$ is dense in $X_{{\alpha}+1}$, we have to guarantee that 
$X_{\alpha}$ is not closed in $X_{{\alpha}+1}$, in particular, $X_{\alpha}$ is not Lindelöf.
(Since in a P-space, Lindelöf subspaces are closed.)
However, in our case it is not automatic in limit steps, because 
the increasing countable union of open non-Lindelöf subspaces can be Lindelöf.

So some extra effort is needed  to  guarantee the non-Lindelöfness in limit steps.

   \vspace{1mm} Assume that $\kappa$ is an uncountable cardinal and $\al$ is a non-zero ordinal. If $X$ is an LLSP space such that $\mbox{ht}(X) = \al$ and $\mbox{wd}(X) = \kappa$, we say that $X$ is a $(\kappa,\alpha)$-{\em LLSP space}.

   \vspace{1mm} Then, we will also transfer the results proved in 
   \cite{bs87} and \cite{ma01}  on thin-tall spaces to the context of locally Lindel\"{o}f P-spaces, showing that Con(ZFC) implies Con(ZFC + ``there is an $(\omega_1,\alpha)$-LLSP space for every ordinal $\alpha < \omega_4$'').

\section{Construction of an LLSP space of width $\om_1$ and height $\om_2$}

   \vspace{2mm} By a {\em decomposition} of a set $A$ of size $\om_1$, we mean a partition of $A$ into subsets of size $\om_1$. In this section we will prove the following result.

\begin{theorem} There is an $(\om_1,\om_2)$-LLSP space. \end{theorem}

\begin{proof}
 We construct an $(\om_1,\om_2)$-LLSP space whose underlying set is $\om_2$. 
For every $\al < \om_2$, we put $I_{\al} = (\om_1 \cdot (\al + 1))\setminus (\om_1 \cdot \al)$, 
and for every ordinal $\xi < \om_1$, we define the ``column'' $N_{\xi} = \{\om_1 \cdot \mu + \xi : \mu < \om_2 \}$. Write ${\xi}\in N_{n({\xi})}$.
Our aim is to construct, by transfinite induction on $\al < \om_2$ an LLSP space $X_{\al}$ satisfying the following:

\vspace{1mm} (1) $X_{\al}$ is an $(\om_1,\al + 1)$-LLSP space such that $I_{\be}(X_{\al}) = I_{\be}$ for every $\be\leq \al$.

\vspace{1mm} (2) For every $\xi < \om_1$, $N_{\xi}\cap X_{\al}$ is a closed discrete subset of $X_{\al}$.

\vspace{1mm} (3) If $\be < \al$ and $x\in X_{\be}$, then a neighbourhood basis of $x$ in $X_{\be}$ is also a neighbourhood basis of $x$ in $X_{\al}$.

\vspace{2mm} For every $\al < \om_2$ and $x\in I_{\al}$, in order to define the required neighbourhood basis ${\mathbb B}_x$ of $x$ in $X_{\al}$, we will also fix a Lindel\"{o}f cone $V_x$ of $x$ in $X_{\al}$ such that the following holds:

\vspace{1mm} (4) $V_x \cap I_{\al} = \{x\}$.

\vspace{1mm} (5) $V_x = \bigcup {\mathbb B}_x$.

\vspace{1mm} (6) There is a club subset $C_x$ of $\om_1$ such that $\om_1 \setminus C_x$ is unbounded in $\om_1$ and $V_x \cap \bigcup\{N_{\nu}: \nu \in C_x \} = \emptyset$.

\vspace{2mm} We define $X_0$ as the set $I_0 = \om_1$ with the discrete topology, and for $x\in I_0$ we put $V_x = \{x\}$ and $C_x = \{y \in \om_1 : y \mbox{ is a limit ordinal } > x \}$. So, assume that $\al > 0$. If $\al = \be + 1$ is a successor ordinal, we put $Z = X_{\be}$. And if $\al$ is a limit ordinal, we define $Z$ as the direct union of
 $\{X_{\be}: \be < \al\}.$ So, the underlying set of the required space $X_{\al}$ is $Z\cup  I_{\al}$. If $x\in Z$, then a basic neighbourhood of $x$ in  $X_{\al}$ is a  neighbourhood of $x$ in  $Z$. Our purpose is to define a neighbourhood basis of each element of $I_{\al}$. Let $\{x_{\nu} : \nu < \om_1 \}$ be an enumeration without repetitions of $Z$. By the induction hypothesis, for every $\xi < \om_1$ there is a club subset $C_{\xi}$ of $\om_1$ such that $\om_1\setminus C_{\xi}$ is unbounded in $\om_1$ and $V_{x_{\xi}} \cap \bigcup \{N_{\nu} : \nu \in C_{\xi} \} = \emptyset$. Let $C = \Delta \{C_{\xi} : \xi < \om_1 \}$, the diagonal intersection of the family  $\{C_{\xi} : \xi < \om_1 \}$. As $V_{x_{\xi}} \cap \bigcup \{N_{\nu} : \nu \in C_{\xi} \} = \emptyset$, by the definition  of $C$, for every $\xi < \om_1$, $V_{x_{\xi}} \cap \bigcup \{N_{\nu} : \nu \in C \}\subset \bigcup \{N_{\nu} : \nu \leq \xi \}$, and clearly $\om_1\setminus C$ is unbounded in $\om_1$. Then, we will define for every element $y\in I_{\al}$ a neighbourhood basis of $y$ from a set $V_y$ in such a way that for some final segment $C'$ of $C$ we will have that  $V_{y} \cap \bigcup \{N_{\nu} : \nu \in C' \} = \emptyset$. We distinguish the following three cases:

 \vspace{5mm}\noindent  {\bf Case 1}. $\al = \be + 1$ is a successor ordinal.

\vspace{3mm}
 For each $\xi < \om_1$ we take a Lindel\"{o}f clopen cone $U_{\xi}$ on some $u_{\xi}$ in $Z$ as follows. We take $U_0\subset V_{x_0}$ as a Lindel\"{o}f clopen cone on $x_0$  such that $(U_0\setminus \{x_0\})\cap N_0 = \emptyset$. Suppose that $\xi > 0$. Let $u_{\xi}$ be the first element $x_{\eta}$ in the enumeration $\{x_{\nu} : \nu < \om_1 \}$ of $Z$ such that $u_{\xi}\not\in \bigcup \{U_{\mu}: \mu < \xi \}$. 
 Since $I_{\beta}\cap \bigcup \{U_{\mu}: \mu < \xi \}\subs \{u_{\mu}:{\mu}<{\xi}\}$, the element
 $u_{\xi}$ is defined. 
 Then, we choose $U_{\xi} \subset V_{x_{\eta}}$ as a Lindel\"{o}f clopen cone on $u_{\xi}$ such that $U_{\xi} \cap \bigcup \{U_{\mu}: \mu < \xi \} = \emptyset$ and $(U_{\xi}\setminus \{u_{\xi}\}) \cap \bigcup \{N_{\nu}: \nu \leq \eta  \} = \emptyset$. So, as  $V_{x_{\eta}} \cap \bigcup \{N_{\nu} : \nu \in C \}\subset \bigcup \{N_{\nu} : \nu \leq \eta \}$, we deduce that $(U_{\xi}\setminus \{u_{\xi}\}) \cap \bigcup \{N_{\nu}: \nu \in C  \} = \emptyset$. And clearly, $\{U_{\xi} : \xi < \om_1 \}$ is a partition of $Z$. Let

$$A = \{\xi \in \om_1 : u_{\xi}\in I_{\be}\cap N_{\rho} \mbox{ for some } \rho \in \om_1\setminus C \}.$$

\noindent 
Since $I_{\beta}\subs \{u_{\xi}:{\xi}<{\omega}_1\}$, we have $|A|={\omega}_1$.
Let $\{A_{\xi}: \xi < \om_1 \}$ be a decomposition of $A$. Fix $\xi < \om_1$. Let $y_{\xi} = \om_1\cdot \al + \xi$. Then, we define

$$V_{y_{\xi}} = \{y_{\xi}\} \cup \bigcup \{U_{\nu} : \nu\in A_{\xi} \}.$$

\noindent Note that since $\bigcup \{U_{\nu} : \nu \in A_{\xi} \}\cap \bigcup \{N_{\nu} : \nu \in C \} = \emptyset$, we infer that $V_{y_{\xi}}\cap \bigcup \{N_{\nu} : \nu \in C \mbox{ and } \nu > \xi \} = \emptyset$. Now, we define a basic neighbourhood of $y_{\xi}$ in $X_{\al}$ as a set of the form

$$\{y_{\xi}\} \cup \bigcup \{U_{\nu} : \nu\in A_{\xi}, \nu\geq\zeta \}$$

\noindent where $\zeta < \om_1$. Then, it is easy to check that conditions $(1)-(6)$ hold.

\vspace{2mm}\noindent  {\bf Case 2}. $\al$ is a limit ordinal of cofinality $\om_1$.

\vspace{2mm} Let $\langle \al_{\nu} : \nu < \om_1 \rangle$ be a strictly increasing sequence of ordinals cofinal in $\al$. For every $\xi < \om_1$, we choose a Lindel\"{o}f clopen cone $U_{\xi}$ on some point $u_{\xi}$ in $Z$ as follows. If $\xi$ is not a limit ordinal, let $u_{\xi}$ be the first element $x_{\eta}$ in the enumeration $\{x_{\nu} : \nu < \om_1 \}$ of $Z$ such that $u_{\xi}\not\in \bigcup \{U_{\mu} : \mu <\xi \}$ and let $U_{\xi}\subset V_{x_{\eta}}$ be a Lindel\"{o}f clopen cone on $u_{\xi}$ such that $U_{\xi}\cap \bigcup \{U_{\mu} : \mu < \xi \} = \emptyset$. Now, assume that $\xi$ is a limit ordinal. Let $\nu < \om_1$ be such that $\al_{\nu} > \mbox{sup} \{\rho(u_{\mu},Z) : \mu < \xi \}$. Then,
we pick $u_{\xi}$ as the first element $x_{\eta}$ in the enumeration $\{x_{\nu} : \nu < \om_1 \}$ of $Z$ such that  $u_{\xi}\in I_{\al_{\nu}}(Z) \cap N_{\delta}$ for some $\delta\in \om_1\setminus C$ with $\delta > \xi$. Note that by the election of $\al_{\nu}$, we have that $u_{\xi}\not\in \bigcup\{U_{\mu} : \mu < \xi \}$. Then, we choose $U_{\xi}\subset V_{x_{\eta}}$ as a Lindel\"{o}f clopen cone on $u_{\xi}$ such that

$$U_{\xi}\cap \bigcup \{U_{\mu} : \mu < \xi \} = \emptyset \mbox{ and }$$

$$(U_{\xi}\setminus \{u_{\xi}\}) \cap \bigcup \{N_{\nu} : \nu \leq \eta \} = \emptyset.$$

\noindent Then since $V_{x_{\eta}} \cap \bigcup \{N_{\nu} : \nu \in C \} \subset \bigcup \{N_{\nu} : \nu \leq \eta \}$ and $\delta\not\in C$, we infer that $U_{\xi}\cap \bigcup \{N_{\nu} : \nu \in C \} = \emptyset$.

\vspace{1mm} Now, let $\{A_{\xi} : \xi < \om_1 \}$ be a decomposition of the set of limit ordinals of $\om_1$. Fix $\xi < \om_1$. Let $y_{\xi} = \om_1 \cdot \al + \xi$. Then, we define

$$V_{y_{\xi}} = \{y_{\xi}\} \cup \bigcup \{U_{\mu} : \mu \in A_{\xi} \}.$$

 \noindent Clearly,
$V_{y_{\xi}} \cap \bigcup \{N_{\nu} : \nu \in C, \nu > \xi \} = \emptyset.$
Now, we define a basic neighbourhood of $y_{\xi}$ in $X_{\al}$ as a set of the form

  $$V_{y_{\xi}} \setminus \bigcup \{U_{\nu} : \nu \in A_{\xi}, \nu < \zeta \}$$

    \noindent where $\zeta < \om_1$.

    \vspace{2mm} Note that the condition that $\delta > \xi$ in the election of $u_{\xi}$ for $\xi$ a limit ordinal is needed to assure that $N_{\xi}\cap X_{\al}$ is a closed discrete subset of $X_{\al}$ for $\xi < \om_1$. So, conditions $(1)-(6)$ hold.

\vspace{3mm}\noindent  {\bf Case 3}. $\al$ is a limit ordinal of cofinality $\om$.

\vspace{2mm} Let $\langle \al_n : n < \om \rangle$ be a strictly increasing sequence of ordinals converging to $\al$. Proceeding by transfinite induction on $\xi < \om_1$, we construct a sequence $\langle u^{\xi}_n : n < \om \rangle$ of points in $Z$ and a sequence $\langle U^{\xi}_n : n < \om \rangle$ such that each $U^{\xi}_n\subset V_{u^{\xi}_n}$ is a Lindel\"{o}f clopen cone on $u^{\xi}_n$ as follows. Fix $\xi < \om_1$, and assume that for $\mu < \xi$ the sequences $\langle u^{\mu}_n : n < \om \rangle$ and $\langle U^{\mu}_n : n < \om \rangle$ have been constructed. Let $C^* = \bigcap \{C_{u^{\mu}_n} : \mu < \xi, n <\om \}$. Note that  $C^*$ is a club subset of $\om_1$, because it is a countable intersection of club subsets of $\om_1$. Now since for every $\mu < \xi$ and $n < \om$, we have that $V_{u^{\mu}_n}\cap \bigcup \{N_{\nu} : \nu \in C_{u^{\mu}_n} \} =\emptyset$, we infer that

$$\bigcup \{V_{u^{\mu}_n}: \mu < \xi, n < \om \} \cap \bigcup\{N_{\nu} : \nu \in C^* \} = \emptyset.$$

\noindent Hence, for every ordinal $\be < \al$,

$$|I_{\be}\setminus \bigcup \{V_{u^{\mu}_n} : \mu < \xi, n <\om \}| = \om_1.$$

\vspace{1mm} Now, we construct the sequences $\langle u^{\xi}_n : n < \om \rangle$ and $\langle U^{\xi}_n : n < \om \rangle$ by induction on $n$. If $n$ is even, let $u^{\xi}_n$ be the first element $x_{\eta}$ in the enumeration $\{x_{\nu} : \nu < \om_1 \}$ of $Z$ such that
$u^{\xi}_n\not\in \bigcup \{U^{\mu}_k : \mu < \xi, k < \om \} \cup \bigcup \{U^{\xi}_k: k < n \}$, and let $U^{\xi}_n\subset V_{x_{\eta}}$ be a Lindel\"{o}f clopen cone on $u^{\xi}_n$  such that

$$U^{\xi}_n \cap (\bigcup \{U^{\mu}_k: \mu < \xi, k < \om \} \cup \bigcup \{U^{\xi}_k:
k < n \}) = \emptyset.$$

Now, suppose that $n$ is odd. Let $k\in \omega$ be such that $\al_k > \mbox{sup}\{\rho(u^{\xi}_m,Z) : m < n \}$. First, we pick $\tilde{u}^{\xi}_n$ as the first element $x_{\eta}$ in the enumeration $\{x_{\nu} : \nu < \om_1 \}$ of $Z$ such that $\tilde{u}^{\xi}_n\in I_{\al_k + 1}(Z)\cap N_{\zeta^*}$ for some $\zeta^*\in C^*$. So, $\tilde{u}^{\xi}_n\not\in \bigcup \{U^{\mu}_m: \mu < \xi, m < \om \} \cup \bigcup \{U^{\xi}_m:
m < n \}$. Now, we choose $\tilde{U}^{\xi}_n\subset V_{x_{\eta}}$ as a  Lindel\"{o}f clopen cone on $\tilde{u}^{\xi}_n$ such that

$$\tilde{U}^{\xi}_n \cap (\bigcup \{U^{\mu}_m: \mu < \xi, m < \om \} \cup \bigcup \{U^{\xi}_m:
m < n \}) = \emptyset.$$

\noindent and

$$(\tilde{U}^{\xi}_n \setminus \{\tilde{u}^{\xi}_n\})\cap \bigcup \{N_{\nu} : \nu \leq \eta \} = \emptyset.$$

\noindent Then as $\tilde{u}^{\xi}_n = x_{\eta}$ and $V_{x_{\eta}} \cap \bigcup\{N_{\nu}: \nu \in C \}\subset \bigcup \{N_{\nu}: \nu \leq \eta \}$, we infer that $(\tilde{U}^{\xi}_n\setminus \{ \tilde{u}^{\xi}_n \})
\cap \bigcup\{N_{\nu}: \nu \in C \} = \emptyset$. However, note that if $\zeta$ is the ordinal such that $\tilde{u}^{\xi}_n \in N_{\zeta}$, it may happen that $\zeta\in C$. Then, we pick $u^{\xi}_n$ as the first element $x_{\rho}$ in the enumeration $\{x_{\nu} : \nu < \om_1 \}$ of $Z$ such that $u^{\xi}_n\in \tilde{U}^{\xi}_n \cap I_{\al_k}(Z)\cap N_{\delta}$ for some $\delta > \xi$. Note that $\delta\not\in C$, because $(\tilde{U}^{\xi}_n\setminus \{ \tilde{u}^{\xi}_n \})\cap \bigcup\{N_{\nu}: \nu \in C \} = \emptyset$. Now, we choose $U^{\xi}_n \subset \tilde{U}^{\xi}_n \cap V_{x_{\rho}}$ as a  Lindel\"{o}f clopen cone on $u^{\xi}_n$ such that

 $$(U^{\xi}_n \setminus \{u^{\xi}_n\})\cap \bigcup \{N_{\nu} : \nu \leq \rho \} = \emptyset.$$

\noindent Hence as $V_{x_{\rho}} \cap \bigcup \{N_{\nu} : \nu \in C \} \subset \bigcup \{N_{\nu} : \nu \leq \rho \}$ and $\delta\not\in C$, we infer that $U^{\xi}_n \cap \bigcup\{N_{\nu}: \nu \in C \} = \emptyset$.

\vspace{1mm} Now, let $\{A_{\xi} : \xi < \om_1 \}$ be a decomposition of $\om_1$. Fix $\xi < \om_1$. Let $y_{\xi} = \om_1 \cdot \al + \xi$. Then, we define

$$V_{y_{\xi}} = \{y_{\xi}\} \cup \bigcup \{U^{\mu}_n : \mu \in A_{\xi}, n \mbox{ odd} \}.$$

\noindent As $\bigcup \{U^{\mu}_n : \mu \in A_{\xi}, n \mbox{ odd} \} \cap \, \bigcup \{N_{\nu} : \nu \in C \} = \emptyset$, we deduce that  $V_{y_{\xi}} \cap \, \bigcup \{N_{\nu} : \nu \in C \mbox{ and } \nu > \xi \} = \emptyset$. Then, we define a basic neighbourhood of $y_{\xi}$ in $X_{\al}$ as a set of the form

$$\{y_{\xi}\} \cup \bigcup \{U^{\mu}_n : \mu\in A_{\xi}, \mu \geq \zeta, n \mbox{ odd} \}$$

\noindent where $\zeta < \om_1$.  Now, it is easy to see that conditions $(1)-(6)$ hold.

 \vspace{2mm} Then, we define the desired space $X$ as the direct union of the spaces $X_{\al}$ for $\al < \om_2$.  
\end{proof}

\vspace{1mm}\noindent {\bf Remark 2.2.} Note that by the construction carried out in the proof of Theorem 2.1, we have that 
\begin{equation}\notag
\text{if $U\subs X$ is Lindel\"{o}f then $\{{\xi}:N_{\xi}\cap U\ne \empt\}\in NS({\omega}_1)$.
}
\end{equation}

\section{A stepping up theorem}

 \vspace{2mm} In this section, for every  cardinal $\la\geq \om_2$  we will construct from an $(\om_1,\la)$-LLSP space satisfying certain additional properties an $(\om_1,\al)$-LLSP space for every ordinal $\al < \la^{+}$. As a consequence of this construction, we will be able to extend Theorem 2.1 from $\om_2$ to any ordinal $\al < \om_3$.  We need some preparation.

 \begin{definitions}
 {\em (a) Assume that $X$ is an LLSP space, $\be + 1 < \mbox{ht}(X)$,  $x\in I_{\be +1}(X)$  and ${\mathbb B}_x$ is a neighbourhood basis for $x$. We say that ${\mathbb B}_x$ is {\em admissible}, if there is a pairwise disjoint family $\{U_{\nu} : \nu < \om_1 \}$ such that for every $\nu < \om_1$, $U_{\nu}$ is a Lindel\"{o}f clopen cone on some point $x_{\nu}\in I_{\be}(X)$  in such a way that ${\mathbb B}_x$ is the collection of sets of the form

 $$\{x\}\cup \bigcup\{U_{\nu} : \nu \geq \xi \},$$

 \noindent  where $\xi < \om_1$.
 Then, we will say that ${\mathbb B}_x$ is the {\em admissible basis for} $x$ {\em given by} $\{U_{\nu} : \nu < \om_1 \}$.

 \vspace{1mm}
 (b) Now, we say that $X$ is an {\em admissible space} if for every $x\in X$ there is a neighbourhood basis ${\mathbb B}_x$ such that for  every successor ordinal $\be + 1 < \mbox{ht}(X)$ the following holds:

  \begin{enumerate}[(1)]

  \item ${\mathbb B}_x$ is an admissible basis  for  every point $x\in I_{\be + 1}(X)$,

  \item if $x,y\in I_{\be + 1}(X) $ with $x\neq y$ and $\rho(x) = \rho(y)$, ${\mathbb B}_x$ is given by $\{U_{\nu} : \nu < \om_1 \}$ and  ${\mathbb B}_y$ is given by $\{U'_{\nu} : \nu < \om_1 \}$, then for every $\nu,\mu < \om_1$ we have $U_{\nu}\cap U'_{\mu} = \emptyset$.

  \end{enumerate}}

 \end{definitions}

\vspace{2mm} Note that the space $X$ constructed in the proof of Theorem 2.1 is admissible.

 \begin{definition} {\em We say that an LLSP space $X$ is {\em good}, if for every ordinal $\al < \mbox{ht}(X)$ and every set $\{U_n : n\in \om \}$ of Lindel\"{o}f clopen cones on points of $X$, the set $I_{\al}(X)\setminus \bigcup\{U_n : n\in \omega \}$ is uncountable.}
 \end{definition}

\vspace{2mm} Note that the space $X$ constructed in the proof of Theorem 2.1 is good.

\vspace{1mm} Assume that $X$ is a good LLSP space. Then, we define the space  $X^*$ as follows. Its underlying set is $X\cup \{z\}$ where $z\not\in X$. If $x\in X$, a basic neighbourhood of $x$ in $X^*$ is a neighbourhood of $x$ in $X$. And a basic neighbourhood  of $z$ in $X^*$ is a set of the form

$$X^*\setminus \bigcup\{U_n : n\in \om \}$$

\noindent where each $U_n$ is a Lindel\"{o}f clopen cone on some point of $X$. Clearly, $X^*$ is a Lindel\"{o}f scattered Hausdorff P-space with $\mbox{ht}(X^*) = \mbox{ht}(X) + 1$.

\begin{theorem} Let $\la \geq \om_2$ be a cardinal. Assume that there is a good $(\om_1,\la)$-LLSP space that is admissible. Then, for  every ordinal $\al < \la^+$ there is a good $(\om_1,\al)$-LLSP space.

\end{theorem}

So, we obtain the following consequence of Theorems 2.1 and 3.3.

\begin{corollary} For every ordinal $\al < \om_3$ there is a good $(\om_1,\al)$-LLSP space.

\end{corollary}

\begin{proof}[Proof of Theorem 3.3]

We may assume that $\la \leq \al < \la^+$. We 
proceed by transfinite induction on $\al$. If $\al = \la$, the case is obvious. Assume that $\al = 
\be + 1$ is a successor ordinal. Let $Y$ be a good $(\om_1,\be)$-LLSP space. For every $\nu < 
\om_1$ let $Y_{\nu}$ be a P-space homeomorphic to $Y^*$ in such a way that $Y_{\nu}\cap Y_{\mu} = 
\emptyset$ for $\nu < \mu < \om_1$. Clearly, the topological sum of the spaces $Y_{\nu}$ ($\nu < 
\om_1$) is a good $(\om_1,\al)$-LLSP space.

\vspace{2mm} Now, assume that $\al > \la$ is a limit ordinal. Let $\theta = \mbox{cf}(\al)$. Note 
that since there is a good admissible $(\om_1,\la)$-LLSP space and $\theta \leq \la$, there is a 
good admissible LLSP space $T$ of width $\om_1$ and height $\theta$.
Let $\{\alpha_{\xi} : \xi < \theta \}$ be a closed strictly increasing sequence of ordinals 
cofinal in $\al$ with $\al_0 = 0$. For every ordinal $\xi < \theta$, we put $J_{\xi} = \{\al_{\xi}
\} \times \om_1$. We may assume that the underlying set of $T$ is $\bigcup \{J_{\xi} : \xi < 
\theta \}$, $I_{\xi}(T) = J_{\xi}$ for every $\xi < \theta$ and $I_{\theta}(T) = \emptyset$. 

Fix a system of neighbourhood bases, $\{\mathbb B_x:x\in T\}$,  which 
witnesses that $T$ is admissible.
Write  $V_s = \bigcup {\mathbb B}_s$ for $s\in T$.

So, writing $$T' = \{ s\in T : \rho(s,T) \mbox{ is a successor  ordinal} \},$$
for each $s\in T$ with  ${\rho}(s,T)={\xi}+1$, there is 
$D_s=\{d^s_{\zeta}:{\zeta}<{\omega}_1\}\in [I_{\xi}(T)]^{{\omega}_1}$
and for each $d\in D_s$ there is a Lindelöf cone  $U_d$ on $d$
such that 
\begin{displaymath}
   \mathbb B_s=\big\{\{s\}\cup \bigcup_{{\eta}\le {\zeta}}U_{d^s_{\zeta}}:{\eta}<{\omega}_1\big\}.
   \end{displaymath}

\vspace{2mm} In order to carry out the desired construction, we will insert an adequate LLSP space 
between $I_{\xi}(T)$ and $I_{\xi + 1}(T)$ for every $\xi < \theta$. If $\xi < \theta$, we define  
$\delta_{\xi} = \mbox{o.t.}(\al_{\xi + 1}\setminus \al_{\xi})$.
We put $y^{\xi + 1}_{\nu} = \langle \al_{\xi + 1},\nu \rangle$ for $\xi < \theta$ and 
$\nu < \om_1$, and we put 
$D^{\xi}_{\nu} = \{x\in T : \rho(x,T) = \xi \mbox{ and } x\in V_{y^{\xi 
+ 1}_{\nu}}\}=D_{y^{{\xi}+1}_{\nu}}$. Since $T$ is admissible, 
$D^{\xi}_{\nu} \cap D^{\xi}_{\mu} = \emptyset$ for $\nu\neq \mu$.

Now, by the induction hypothesis, for every point $y = y^{\xi + 1}_{\nu}$ where $\xi < \theta$ and 
$\nu < \om_1$ there is a Lindel\"{o}f scattered Hausdorff P-space $Z_y$ of height $\de_{\xi} + 1$ 
such that $I_0(Z_y) = D^{\xi}_{\nu}$, 
$|I_{\nu}(Z_y)|={\omega}_1$ for ${\nu}<{\delta}_{\xi},$
$I_{\de_{\xi}}(Z_y) = \{y\}$ and $Z_y \cap T = D^{\xi}_{\nu} 
\cup \{y\}$. 
Also, we assume that $Z_{y^{\xi + 1}_{\nu}} \cap Z_{y^{\xi + 1}_{\mu}} = \emptyset$ for 
$\nu\neq\mu$ and 
$(Z_{y^{\xi + 1}_{\nu}}\setminus \{y^{\xi + 1}_{\nu}\}) \cap (Z_{y^{\eta + 1}_{\mu}}\setminus \{y^
{\eta + 1}_{\mu}\}) = \emptyset$ 
for $\xi\neq \eta$ and $\nu,\mu < \om_1$.

\vspace{2mm} Now, our aim is to define the desired $(\om_1,\al)$-LLSP space $Z$. Its underlying set is

$$Z=T \cup \bigcup\{Z_y : y\in T'\}.$$

\noindent If $V$ is a Lindel\"{o}f clopen cone on a point $z\in T$, we define

$$V^* = V \cup \bigcup\{ (Z_y\setminus T) : y\in V\cap T' \}.$$
Observe that if $y\in V\cap T'$, then $Z_y\setm V^*= D_y\setm V$
and $D_y\setm V$ is countable because $T$ is admissible. So $Z_y\cap V^*$ is open in $Z_y$
because $Z_y$ is a P-space.

 Now, assume that $x\in Z_s$ for some $s\in T'$. 
Then, if $U$ is a Lindel\"{o}f 
clopen cone on $x$ in $Z_s$, we define

$$U^{\sim} = U \cup \bigcup \{(U_y)^* : y\in D_s\cap U \}.$$

\vspace{2mm} Note that for every $s\in T'$ we have $(V_s)^* = (Z_s)^{\sim}$.

\vspace{2mm}
After that preparation we can define the bases of the points of $Z$.
Suppose that $x \in Z  =  T \cup \bigcup\{Z_y : y\in T'\}$.

If $x\in T\setm T'$, 
then let 
\begin{displaymath}
\mathbb B^Z_x=\{V^*:\text{$V$ is a Lindel\"{o}f clopen cone on $x$ in $T$}\}.
\end{displaymath}

If $x\in (Z\setm T)\cup T'$,
then pick first the unique $s\in T'$ such that $x\in Z_s\setm I_0(Z_s)$, and  
let 
\begin{displaymath}
   \mathbb B^Z_x=\{U^\sim :\text{$U$ is a Lindel\"{o}f clopen cone on $x$ in $Z_s$}\}.
   \end{displaymath}

\noindent{\bf Claim 1.}
{\em $\{\mathbb B^Z_x:x\in Z\}$ is a system of neighbourhood bases of a topology ${\tau}_Z$.}

\begin{proof}
Assume that $y\in W\in \mathbb B^Z_x$.   
We should show that  $\mathbb B^Z_y\cap \mathcal P(W)\ne \empt$.

Assume first that  $x\in T\setm T'$, and so $W=V^*$ for some Lindelöf clopen clone $V$ on $x$ in $T$.

If $y\in T\setm T'$, then $y\in V$ and so   $S\subs V$ for some Lindelöf clopen clone $
S$ on $y$ in $T$.
Thus $y\in S^*\subs V^*$ and $S^*\in \mathbb B^Z_y$.

If $y\in (Z\setm T)\cup T'$ then pick first the unique $s\in T'$ such that 
$y\in Z_s\setm I_0(Z_s)$. 
Then $s\in V$ because otherwise $y\in V^*$ is not possible. 
So as we observed, $V^*\cap Z_s$ is open in $Z_s$. 
So let $S$ be a Lindel\"{o}f clopen cone on $y$ in $Z_s$ with $S\subs V^*\cap Z_s$. 
Then $y\in S^\sim\subs V^*$ and $S^\sim\in \mathbb B^Z_y$.

\smallskip
Assume now  that $x\in (Z\setm T)\cup T'$, then pick first the unique $s\in T'$ such that $x\in Z_s\setm I_0(Z_s)$. Then $W=U^\sim $ for some Lindelöf clopen cone $U$ on $x$ in $Z_s$.

If $y\in Z_s\setm I_0(Z_s)$, then  $S\subs U$ for some Lindelöf clopen clone $S$ on $y$ in $Z_s$,
and so $S^\sim \in \mathbb B^Z_y$ and $S^\sim\subs U^\sim$.

If $y\notin Z_s \setm I_0(Z_s)$, then $y\in (U_d)^*$ for some $d\in I_0(Z_s)\cap U$, and so 
there is $S\in \mathbb B^Z_y$ with $S\subs (U_d)^*$ using what we proved so far. 
Thus $S\subs U^\sim$ as well. 
\end{proof}

\noindent{\bf Claim 2.}
{\em ${\tau}_Z$ is Hausdorff.}

\begin{proof}
Assume that $\{x,y\}\in {[Z]}^{2}$. Let $s$ and $t$ be elements of $T$ such that 
$x\in Z_s\setm I_0(Z_s)$ if $x\notin T\setm T'$ and $s=x$ otherwise, and  
$y\in Z_t\setm I_0(Z_t)$ if $y\notin T\setm T'$ and $t=y$ otherwise.

If $s\neq t$, consider disjoint Lindel\"{o}f clopen cones $U$ and $V$ on $s$ and $t$  in $T$ respectively. Note that if $w\in U \cap T'$, then $Z_w\setminus T \subset U^*$ because $w\in U$, but $(Z_w\setminus T)\cap V^* = \emptyset$ because $w\not\in V$, and analogously if $w\in V \cap T'$ then $Z_w\setminus T \subset V^*$ but $(Z_w\setminus T)\cap U^* = \emptyset$. So, $U^*$ and $V^*$ are disjoint open sets containing $x$ and $y$ respectively.

If $s=t$, then there are disjoint cones in $Z_s$ on $x$ and $y$,  $U$ and $V$, respectively.
Then $U^*$ and $V^*$ are disjoint open sets containing $x$ and $y$, respectively. 
\end{proof}

It is trivial from the definition that $Z$ is a $P$-space because 
$T$ is a P-space and the  $Z_s$ are P-spaces.

By transfinite induction on ${\delta}<{\alpha}$  
 it is easy to check that  
 \begin{displaymath}
 {I_{\delta}(Z)}=\left\{\begin{array}{ll}
 {J_{\xi}}&\text{if ${\delta}={\alpha}_{\xi}$},\\\\
 {\bigcup\{I_{\eta}(Z_s\}:s\in I_{{\xi}+1}(T)\}}&\text{if ${{\alpha}_{\xi}<\delta}={\alpha}_{\xi}+{\eta}<{\alpha}_{\xi+1}$,}\\
 \end{array}\right.
 \end{displaymath}
 so $Z$ is scattered with  height $\al$ and width ${\omega}_1$.

\smallskip

\smallskip

\noindent{\bf Claim 3.} {\em $Z$ is locally Lindelöf.}

\begin{proof}
   \vspace{1mm} Note that if $x\in T\setminus T'$ and $U^*\in \mathbb B^Z_x$, then for every $V^*\in \mathbb B^Z_x$ with $V^*\subset U^*$ we have that $U^*\setminus V^* = \bigcup \{W^*_n : n \in \nu \}$ where $\nu \leq \omega$ in such a way that each $W_n$  is a Lindel\"{o}f clopen cone on some point
   $v_n\in T \cap U$ in $T$ with $\rho(v_n, T) < \rho(x,T)$.
   
   \vspace{1mm} Also, if $x\in T' \cup (Z\setminus T)$ and $U^{\sim}\in \mathbb B^Z_x$ then for every $V^{\sim}\in \mathbb B^Z_x$ with $V^{\sim}\subset U^{\sim}$, if $s$ is the element of $T'$ with $x\in Z_s\setminus I_0(Z_s)$, we have that $U^{\sim}\setminus V^{\sim} = \bigcup \{U'_n : n\in \nu \}$ where $\nu\leq \omega$ in such a way that for every $n\in \nu$, either $U'_n = U^{\sim}_n$ where $U_n$ is a Lindel\"{o}f clopen cone on some point $u_n \in Z_s\cap U$ in $Z_s$ with $0 < \rho(u_n,Z_s) < \rho(x,Z_s)$ or $U'_n = U^*_n$ where $U_n$ is a Lindel\"{o}f clopen cone on some point $u_n\in D_s\cap U$ in $T$.
   
   \vspace{1mm} Now, proceeding by transfinite induction on $\rho(x,Z)$, we can verify that if $x\in T\setminus T'$ and $U$ is a Lindel\"{o}f clopen cone on $x$ in $T$, then $U^*$ is a Lindel\"{o}f clopen cone on $x$ in $Z$, and that if $x\in Z_s\setminus I_0(Z_s)$ for some $s\in T'$ and $U$ is a Lindel\"{o}f clopen cone on $x$ in $Z_s$, then $U^{\sim}$ is a Lindel\"{o}f clopen cone on  $x$ in $Z$. Therefore, $Z$ is locally Lindel\"{o}f.
   \end{proof}

\noindent{\bf Claim 4.} {\em $Z$ is good.}

\begin{proof}
   Let ${\delta}<{\alpha}=ht(Z)$
and let $\{W_n : n\in \om \}$ be a family of Lindelöf cones in $Z$. Since every $W_n$ is covered by countably many
Lindelöf cones from the basis, we can assume that $W_n\in \mathbb B^Z_{x_n}$ for some $x_n\in Z$
for each $n\in {\omega}$.
For each $n$ pick $y_n\in T$ such that $y_n=x_n$ if $x_n\in T$ and $x_n\in Z_{y_n}$ otherwise.

 Then $W_n\subs W'_n$ for some $W'_n\in \mathbb B^Z_{y_n}$, so we can assume that 
 $\{x_n:n\in {\omega}\}\subs T$. 
We can also assume that if $x_n\in T'$, then $W_n$ is as large as possible, i.e.   
$W_n=Z_{x_n}^\sim=(V_{x_n})^*$.

If $x_n\in T\setm T'$, then $W_n=S_n^*$ for some Lindelöf cone $S_n$ on $x_n$ in $T.$

If ${\delta}={\alpha}_{\xi}$ for some ${\xi}$, then  
$I_{\delta}(Z)\cap W_n=I_{\delta}(Z)\cap V_{x_n}$ if $x_n\in T'$
and $I_{\delta}(Z)\cap W_n=I_{\delta}(Z)\cap S_n$ if $x_n\in T\setm T'$.

So $I_{\delta}(Z)\setm \bigcup_{n\in {\omega}}W_n$ is uncountable because 
$T$ is good. 
 
Assume that  ${\alpha}_{\xi}< {\delta}<{\alpha}_{{\xi}+1}$
and let ${\delta}={\alpha}_{\xi}+{\eta}$.

Pick $s\in I_{{\alpha}_{{\xi}+1}}(Z)\setm \bigcup_{n\in {\omega}}W_n$. Then $Z_s\setm \bigcup_{n\in {\omega}}W_n\supset I_{\eta}(Z_s)$, and so  $I_{\delta}(Z)\setm \bigcup_{n\in {\omega}}W_n\supset I_{\eta}(Z_s)$, and hence $I_{\delta}(Z)\setm \bigcup_{n\in {\omega}}W_n$ is uncountable.

\end{proof}

 Thus, the space $Z$ is as required. 
\end{proof}

\section{Cardinal sequences of length $<\om_4$}

\vspace{2mm} In this section, we will show the following result.

\begin{theorem}\label{tm:o4} If V=L, then there is a cardinal-preserving partial order ${\mathbb P}$ such that in $V^{\mathbb P}$ there is an
  $(\om_1,\al)$-LLSP space for every ordinal $\al < \om_4$.

\end{theorem}

\vspace{2mm} If $S = \bigcup \{ \{\al\} \times A_{\al} : \al < \eta \}$ where $\eta$ is a non-zero ordinal and each $A_{\al}$ is a non-empty set of ordinals, then for every $s = \langle \al,\xi \rangle\in S$ we write $\pi(s) = \al$ and $\zeta(s) = \xi$.

\vspace{1mm} The following notion is a refinement of a notion used implicitly in \cite{bs87}.

\begin{definition} {\em  We say that ${\mathbb S} = \langle S,\preceq, i \rangle$ is an {\em LLSP poset}, if the following conditions hold:

\begin{enumerate}[(P1)]

  \item $\langle S,\preceq \rangle$ is a partial order with $S= \bigcup \{S_{\al} : \al < \eta \}$ for some non-zero ordinal $\eta$ such that each $S_{\al} = \{\al\} \times A_{\al}$ where $A_{\al}$ is a non-empty set of ordinals.

 \item If $s \prec t$ then $\pi(s) < \pi(t)$.

  \item If $\al < \be < \eta$ and $t\in S_{\be}$, then $\{s\in S_{\al} : s \prec t \}$ is uncountable.


 \item  
 If $\ga < \eta$ with $\mbox{cf}(\ga) = \om$,  $t\in S_{\gamma}$ and $\langle t_n : n \in \om \rangle$ is a sequence of elements of $S$ such that $t_n \prec t$ for every $n\in \om$, then for every ordinal $\be < \ga$ the set $\{s\in S_{\be} : s \prec t \mbox{ and } s\not\preceq t_n \mbox{ for } n\in \om \}$ is uncountable.

 \item  $i : [S]^2 \rightarrow [S]^{\leq {\om}}$ such that for every $\{s,t\}\in [S]^2$ the following holds:

     \begin{enumerate}[(a)]
     \item If $v\in i\{s,t\}$ then $v\preceq s,t$.
     \item If $u\preceq s,t$, then there is $v\in i\{s,t\}$ such that $u\preceq v$.
         \end{enumerate}
          \end{enumerate}}
          \end{definition}

          \noindent If there is an uncountable cardinal $\lambda$ such that $|S_{\al}| = \lambda$ for $\al < \eta$, we will say that $\langle S,\preceq, i \rangle$ is a $(\lambda,\eta)$-{\em LLSP poset}.

          \vspace{2mm} If ${\mathbb S} = \langle S,\preceq, i \rangle$ is an LLSP poset 
          with $S= \bigcup \{S_{\al} : \al < \eta \}$, 
          we define its {\em associated LLSP space} $X = X({\mathbb S})$ as follows. 
          The underlying set of $X({\mathbb S})$ is $S$. If $x\in S$ we write $U(x) = \{y\in S: y \preceq x \}$. Then, for every $x\in S$ we define a basic neighbourhood of $x$ in $X$ as a set of the form $U(x)\setminus \bigcup\{U(x_n) : n\in \om \}$ where each $x_n \prec x$. It is easy to check that $X$ is a locally Lindel\"{o}f scattered Hausdorff  P-space (see \cite{ba2007} for a parallel proof). 
          And by conditions $(P3)$ and  $(P4)$ in Definition 4.2, we infer that $\mbox{ht}(X) = \eta$ and 
          $I_{\al}(X) = S_{\al}$ for every $\al < \eta$.

          \vspace{2mm} In order to prove Theorem 4.1, first we will 
          construct an $(\om_1,\om_3)$-LLSP space $X$ in a generic extension by means of 
           an $\om_1$-closed $\om_2$-c.c. forcing, by using an argument similar to the one 
           given by Baumgartner and Shelah in \cite{bs87}.

          \vspace{2mm} Recall that a function $F:[\om_3]^2 \rightarrow [\om_3]^{\leq\om_1}$ has {\em property $\Delta$}, if $F\{\al,\be\}\subset \mbox{min} \{\al,\be\}$ for every $\{\al,\be\}\in [\om_3]^2$ and for every set $D$ of countable subsets of $\om_3$ with $|D| = \om_2$ there are $a,b\in D$ with $a\neq b$ such that for every $\al \in a\setminus b$, $\be \in b\setminus a$ and $\tau \in a\cap b$ the following holds:

          \begin{enumerate}[(a)]
          \item if $\tau < \al,\be$ then $\tau\in F\{\al,\be\}$,
          \item if $\tau < \beta$ then $F\{\al,\tau\} \subset F\{\al,\be\}$,
          \item if $\tau < \al$ then $F\{\tau,\be\} \subset F\{\al,\be\}$.
          \end{enumerate}

          \vspace{1mm} By a result due to Velickovic, it is known that $\square_{\om_2}$ implies the existence of a function $F:[\om_3]^2 \rightarrow [\om_3]^{\leq\om_1}$ satisfying property $\Delta$ (see \cite[Chapter 7 and Lemma 7.4.9.]{to},  for a proof ).

\begin{proof}[Proof of Theorem \ref{tm:o4}]
Let $F:[\om_3]^2 \rightarrow [\om_3]^{\leq\om_1}$ be a function with property $\Delta$. First, we construct by forcing an $(\om_1,\om_3)$-LLSP poset. Let $S= \bigcup \{S_{\al} : \al < \om_3 \}$ where $S_{\al} = \{\al\}\times \om_1$ for each $\al < \om_3$. $S$ will be the underlying set of the required poset. We define $P$ as the set of all $p = \langle x_p,\preceq_p,i_p\rangle$ satisfying the following conditions:

\begin{enumerate}[(1)]
\item $x_p$ is a countable subset of $S$.
\item $\preceq_p$ is a partial order on $x_p$ such that:

\begin{enumerate}[(a)]
\item if $s\prec_p t$ then $\pi(s) < \pi(t)$,
\item if $s\prec_p t$ and $\pi(t)$ is a successor ordinal $\be + 1$, then there is $v\in S_{\be}$ such that $s\preceq_p v  \prec_p t$.
    \end{enumerate}

    \item $i_p : [x_p]^2 \rightarrow [x_p]^{\leq \om }$ satisfying the following conditions:

        \begin{enumerate}[(a)]
        \item if $s\prec_p t$ then $i_p\{s,t\} = \{s\}$,
        \item if $s\not\preceq_p t$ and $\pi(s) < \pi(t)$, then $i_p\{s,t\}\subset \bigcup \{S_{\al} : \al \in F\{\pi(s),\pi(t)\}\}$,
            \item if $s,t\in x_p$ with $s\neq t$ and $\pi(s) = \pi(t)$ then $i_p\{s,t\} = \emptyset$,
                \item $v\preceq_p s,t$ for all $v\in i_p\{s,t\}$,
                \item for every $u\preceq_p s,t$ there is $v\in i_p\{s,t\}$ such that $u\preceq_p v$.
                    \end{enumerate}
                    \end{enumerate}

                    \vspace{1mm} If $p,q\in P$, we write $p\leq q$ iff $x_q \subset x_p$, $\preceq_p \upharpoonright x_q = \preceq_q$ and $i_p\upharpoonright [x_q]^2 = i_q$. We put ${\mathbb P} = \langle P,\leq \rangle$.

                    \vspace{2mm} Clearly, ${\mathbb P}$ is $\om_1$-closed. And since the function $F$ has property $\Delta$, it is easy to check that   ${\mathbb P}$ has the $\om_2$-c.c., and so  ${\mathbb P}$ preserves cardinals.

                    \vspace{2mm} Now, let $G$ be a ${\mathbb P}$-generic filter. 
                    We write $\preceq = \bigcup\{\preceq_p : p\in G \}$ and $i = \bigcup \{i_p : p\in G\}$. It is easy to see that 
                    $S = \bigcup\{x_p : p\in G \}$ and $\preceq$ is a partial order on $S$. 
                    Then, we have that $\langle S,\preceq, i \rangle$ is an $(\om_1,\om_3)$-LLSP poset. For this, note that conditions $(P1),(P2),
                    (P5)$ in Definition 4.2 are obvious, and condition $(P3)$ follows from a basic density argument. So, we verify condition $(P4)
                    $. For every $t\in S$ such that $\gamma = \pi(t)$ has cofinality $\om$, 
                    for every sequence $\langle t_n : n\in \om \rangle$ of elements of 
                    $S$, 
                    for every ordinal $\be < \ga$ and for every ordinal $\xi < 
                    \om_1$ let

                    \vspace{2mm}
                    $D_{t,\{t_n : n\in \om \},\be,\xi} = \{ q\in P : \{t\}\cup \{t_n:n\in \om \} \subset x_q \mbox{ and either } (t_n\not\prec_q t \mbox{ for some } n\in \om) \mbox{ or } (t_n \prec_q t \mbox{ for every } n\in \om \mbox{ and there is } y\in S_{\be}\,\cap\, x_q \mbox{ with }$ $\zeta(y) > \xi \mbox{ such that } y \prec_q t \mbox{ and } y\not\preceq_q t_n \mbox{ for every } n\in \om ) \}.$

                    \vspace{2mm} Since ${\mathbb P}$ is $\om_1$-closed, we have that $D_{t,\{t_n : n\in \om \},\be,\xi}\in V$. Then, consider $p = \langle x_p,\preceq_p,i_p\rangle\in P$. We define a $q\in D_{t,\{t_n : n\in \om \},\be,\xi}$ such that $q\leq p$.  Without loss of generality, we may assume that $t\in x_p$. We distinguish the following cases.

                    \vspace{2mm}\noindent  {\bf Case 1}. $t_n\not\in x_p$ for some $n\in \om$.

\vspace{1mm} We define $q =\langle x_q,\preceq_q,i_q \rangle$ as follows:

\vspace{1mm} (a) $x_q = x_p \cup \{t_n : n\in \om \}$,

\vspace{1mm} (b) $\prec_q = \prec_p$,

\vspace{1mm} (c) $i_q\{x,y\} = i_p\{x,y\}$ if $\{x,y\}\in [x_p]^2$, $i_q\{x,y\} = \emptyset$ otherwise.

 \vspace{2mm}\noindent  {\bf Case 2}. $t_n\in x_p$ for every $n\in \om$.

 \vspace{1mm} If $t_n \not\prec_p t$ for some $n\in \omega$, we put $q = p$. So, assume that $t_n \prec_p t$ for all $n\in \om$. Let $u\in S_{\be}\setminus x_p$ be such that $\zeta(u) > \xi$. We define $q =\langle x_q,\preceq_q,i_q \rangle$ as follows:

\vspace{1mm} (a) $x_q = x_p \cup \{u\}$,

\vspace{1mm} (b) $\prec_q = \prec_p \cup \{\langle u,v \rangle : t \preceq_p v \}$,

\vspace{1mm} (c) $i_q\{x,y\} = i_p\{x,y\}$ if $\{x,y\}\in [x_p]^2$, $i_q\{x,y\} = \{x\}$ if $x \prec_q y$, $i_q\{x,y\} = \{y\}$ if $y \prec_q x$, $i_q\{x,y\} = \emptyset$ otherwise.

\vspace{2mm} So, $D_{t,\{t_n : n\in \om \},\be,\xi}$ is dense in ${\mathbb P}$, and hence condition $(P4)$ holds.  Let $X = X(\langle S,\preceq,i \rangle )$. For every $x\in S$, we write $U(x) = \{y\in S : y\preceq x \}$. By conditions $(2)(b)$ and $(3)(c)$ in the definition of P, we see that if $x\in S_{\be + 1}$ for some $\be < \om_3$, then $x$ has an admissible basis in $X$ given by $\{U(y) : y \prec x, \pi(y) = \be \}$. Thus, $X$ is an admissible space. And clearly, $X$ is good. So, by  Theorem 3.3, we can construct from the space $X$ an $(\om_1,\al)$-LLSP space for every ordinal $\om_3\leq \al < \om_4$. 
\end{proof}

Now, assume that $\ka$ is an uncountable regular cardinal. Recall that a topological space $X$ is a $P_{\ka}$-{\em space}, if the intersection of any family of less than $\ka$ open subsets of $X$ is open in $X$. And we say that $X$ is $\kappa$-{\em compact}, if every open cover of $X$ has  a subcover of size less $< \kappa$.
 By an $SP_{\kappa}$ {\em space} we mean a scattered Hausdorff $P_{\ka}$-space. Then, we want to remark that by using arguments that are parallel to the ones  given in the proofs of the above theorems, we can show the following more general results:

\vspace{1mm} (1) For every uncountable regular cardinal $\ka$ and every ordinal $\al < \ka^{++}$, there is a locally $\ka$-compact  $SP_{\kappa}$ space $X$ such that $\mbox{ht}(X) = \al$ and $\mbox{wd}(X) = \kappa$.

\vspace{1mm} (2) If V=L  and $\ka$ is an uncountable regular cardinal, then there is a cardinal-preserving partial order ${\mathbb P}$  such that in $V^{\mathbb P}$ we have that for every ordinal $\al < \ka^{+++}$ there is a locally $\ka$-compact $SP_{\kappa}$ space $X$ such that $\mbox{ht}(X) = \al$ and $\mbox{wd}(X) = \kappa$.

\end{document}